\documentclass[10pt,leqno]{amsart}
\usepackage{graphicx}
\usepackage{indentfirst,csquotes}

\usepackage{amssymb,amsthm,amsmath}
\usepackage{xcolor,paralist,hyperref,titlesec,fancyhdr,etoolbox}
\def\per{{{\operatorname{per}}}}

\newtheorem*{Embedding Theorem}{Embedding Theorem}
\newtheorem*{theorem 1}{Theorem 1}
\newtheorem*{corollary 1}{Corollary 1}
\newtheorem*{corollary 2}{Corollary 2}
\usepackage{tikz}
\usetikzlibrary{arrows}

\titleformat{\section}[display]{\normalfont\huge\bfseries\centering}{\centering\chaptertitlename\thechapter}{10pt}{\Large}
\titlespacing*{\section}{0pt}{0ex}{0ex}

\hypersetup{ colorlinks=true, linkcolor=black, filecolor=black, urlcolor=black }

\begin{document}
\title{On the subsystems of certain sofic shifts}
\author[Wolfgang Krieger]{Wolfgang Krieger}
\date{\today}
\maketitle

\let\thefootnote\relax

\begin{abstract}
For an aperiodic subshift of finite type $Y$ and for a subshift $X$ with topological entropy less than the topological entropy of $Y$, a theorem is proved in Krieger: On the subsystems of topological Markov chains, Ergodic Theory \& dynamical systems 1982 $\bold{2}$, 195-202, that says that the necessary condition on the periodic points of $X$ and $Y$ for the existence of an embedding of $X$ into $Y$ is also sufficient for the existence of an embedding of $X$ into $Y$. In this note we point out that this theorem extends to  certain classes of sofic shifts as target shifts.
\end{abstract} 

Let $\Sigma$ be a finite alphabet, and let $S_\Sigma$ denote the left shift on 
$\Sigma^\Bbb Z$ Given an $S_\Sigma$-invariant subset $X$ of $\Sigma^\Bbb Z$ the restriction $S_\Sigma$ of $S_\Sigma$ to $X$ is called a subshift.
 For basic notions of the theory of subshifts we refer to \cite {Ki} or \cite {LM}.
Given a subshift $X$, we denote its topological entropy of $X$ by $h(X)$. We denote the set of periodic points of $X$ by $P(X)$, and we denote the least period of a periodic point  of $X$ by $\per (p)$.

In \cite {Kr1}
the following theorem was proved.
  
\begin{Embedding Theorem} 
Let $Y$ be an aperiodic subshift of finite type and let  $X$ be a subshift such that $h(X) < h(Y)$. Then there exists an embedding of $X$ into $Y$ if and only if
$$
\# \{p \in P(X): \per(p) = n \}\leq \# \{p \in P(Y) : \per(p) = n\},
 \quad \quad \quad  n \in \Bbb N.
$$
\end{Embedding Theorem}

In \cite {B} there is an example that shows that the Embedding Theorem does not extend to all mixing sofic shifts as target shifts. In this note we identify two families of sofic shifts to which the Embedding Theorem extends.

We introduce notation and terminology. We denote a finite directed graph  with
vertex set ${\frak V}$ and edge set 
${\mathcal E}$ by $G(\frak V, \mathcal E)$.
A $\Sigma$-labelled directed graph  $G(\frak V, \mathcal E)$ is called 1-right resolving if for every vertex $V \in \frak V$  and every label $\sigma\in \Sigma$ there is at most one edge that leaves $V$ and that carrries the label $\sigma$.
The edge shift of the 1-right resolving $\Sigma$-labelled directed graph 
$G(\frak V, \mathcal E)$ is the set of bi-infinite directed paths in the graph that is equipped with the left shift. 

Given a subshift $X \subset \Sigma^\Bbb Z$, and $I^-, I^+ \in \Bbb Z, I^-\leq I^+$,
we call
$(\sigma_i)_{I^-\leq i \leq I^+}\in \Sigma^{[I^-, I^+]}$ a block
and we use the notation
$$
x_{[I^-, I^+]} = (x_i)_{I^- \leq i \leq I^+},\quad \quad \quad    x \in X,
$$
and the notation
$$
X_{[I^-, I^+]}   = \{x_{[I^-, I^+]}  : x \in A\}.
$$
We use similar notations in the case that indices range in left infinite or right infinite intervals. We use the same notation for a block in $ X_{[I^-, I^+]}$ and for the word that the block carries. 

The language of admissible words of a subshift $X$ we denote by $\mathcal L(X)$.
We set
$$
\Gamma^+_\infty(a) = \{ x^- \in X_{(j, \infty)}: a x^{-}
 \in X_{[i, \infty)}\}, \quad 
a \in  x_{[i,j]}, \quad i,j \in \Bbb Z, i \leq j.
$$
$$
\Gamma^+(a) = \{ b \in \mathcal L(X): ab \in \mathcal L (X)\}, \quad 
a \in  \mathcal L (X).
$$
The notations $\Gamma^{-}, \Gamma^+$ have the symmetric meaning. A word 
$a \in  \mathcal L (X)$ is said to be synchronizing 
if $b \in \Gamma^-(a)$, and $c \in\Gamma^+(a)$ together imply that 
$bac \in  \mathcal L (X)$.

The right Fischer cover \cite{F} of a sofic shift $Y \subset  \Sigma^\Bbb Z$, for which we use the notation
 $\widehat {Y}$, is the edge shift of a $\Sigma$-labelled 1-right resolving directed graph 
 $G(\frak V_Y. \mathcal E_Y)$,  that has as its vertex set 
 $\frak V(Y)$ the set $\Gamma^+(c), c$ a synchronizing words of $Y$. 
 A vertex $V \in\frak V_Y$ accepts the symbol $\sigma \in \Sigma$ as input precisely
 if there is a right-infinite sequence in $V$ that starts with $\sigma$, and the target vertex of the edge in $\mathcal E_Y$ that leaves $V$ and that carries the label $\sigma$ is the set of right-infinite sequences $(\sigma_k)_{k>1}$ that are obtained from a right-infinite sequence 
 $(\sigma_k)_{k\in \Bbb N}$
 such that $\vartheta_1 = \sigma$ by removing its initial symbol $\sigma$.
We denote by $\eta_Y$ the projection of $\widehat{Y}$ onto $Y$ that assigns to a bi-infinite path in $G(\frak V_Y. \mathcal E_Y)$ its label sequence. 
 The sofic shift $Y$ is said to be almost Markov if the graph 
 $G(\frak V_Y. \mathcal E_Y)$ is left resolving. 
 
In \cite{M}
MacDonald has proved an extension of the Embedding Theorem. This extension implies a  theorem for the Fischer cover that we state as Theorem 1

\begin{theorem 1}
Let $Y$ be an aperiodic sofic shift.
Then there exist an embedding $\varphi$ of $X$ into $\widehat{Y}$ such that the restriction of $\eta_{\widehat{Y}}$ to $\varphi (X)$ embeds $\varphi(X)$ into $Y$, if and only if
$$
\# \{p \in P(X): \per(p) = n \}\leq \# \{\eta_{\widehat{Y}}(\widehat{p}): 
\widehat{p} \in P(\widehat{Y}), \per(\widehat{p}) = n\},
 \quad \quad \quad  n \in \Bbb N.
$$
\end{theorem 1}

In \cite{MMTW} there is an extension of MacDonald's result (and therefore also of Theorem 1) to the case of irreducible sofic target shifts. There are two corollaries:
\begin{corollary 1}
The Embedding Theorem extends to irreducible sofic shifts $Y$ as target shifts that have a left (right) Fisher cover in which for all $n\in \Bbb N$ every periodic point of least period $n$ of $Y$ has under $\eta_{\widehat{Y}}$ an inverse image  in $\widehat{Y}$ of least period $n$ .
\end{corollary 1}
\begin{corollary 2}
The Embedding Theorem extends to irreducible sofic shifts as target shifts that have a left (right) Fisher cover that coincides with the left (right) Krieger cover.
\end{corollary 2}
Note that in the corollaries the hypothesis is invariant under topological conjugacy.
The following graph1 is an example of an almost Markov shift that satisfies the hypothesis of Corollary 1:
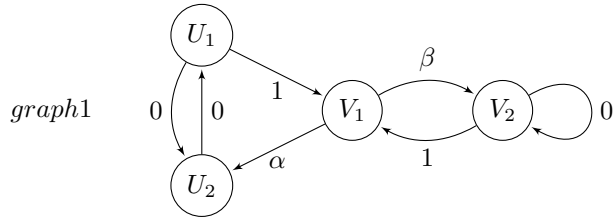
\begin{center}
\begin{tikzpicture}[shorten >=1pt]

\tikzset{vertex/.style = {shape=circle,draw,minimum size=2em}}
\tikzset{edge/.style = {->,> = latex'}}

\node at (0,0) {$graph1$};

\node[vertex] (u1) at  (2,1) {$U_1$};
\node[vertex] (u2) at  (2,-1) {$U_2$};
\node[vertex] (v1) at  (4,0) {$V_1$};
\node[vertex] (v2) at  (6,0) {$V_2$};

\draw[edge] (v1) to[bend left] node[midway, above] {$\beta$} (v2);
 
\draw[edge] (v1) to node[midway, below]{$\alpha$} (u2);

\draw[edge] (u1) to node[midway, below]{$1$} (v1);

\draw[edge] (u2) to node[midway, right]{$0$} (u1);

\draw[edge] (u1) to[bend right] node[midway, left] {$0$} (u2);

\draw[edge] (v2) to[loop, out=30, in=-30, looseness=8] node[right] {$0$} (v2);

\draw[edge] (v2) to[bend left] node[below] {$1$} (v1);

\end{tikzpicture}
\end{center}

Given a subshift $X$ and $l, K \in \Bbb N, K \geq 2l,$ we say that a word 
$$
c = (c_i)_{1 \leq i \leq K} \in \mathcal L(X) 
$$
is periodic of period $l$ if
$$
c_i=c_{i+ l},        \quad \quad \quad (1 \leq i \leq K - l).
$$
In the case that for a point $c\in \mathcal L(X)$ of period $l$ there exists a periodic point $p$ of $X$ such that
$$
p_{[0, l)} = c_{[0, l)},
$$ 
this periodic point is unique.

To emphasize the simplicity of the coding argument we give first the proof of the Embedding Theorem. In conclusion we then indicate the addition to the proof that are necessary to obtain Theorem 1.

We restate and prove the Embedding theorem.
\begin{Embedding Theorem}
Let $Y\subset \Sigma^\Bbb N$ be an aperiodic subshift of finite type and let  $X $ be a subshift such that $\Delta = h(Y) - h(X)> 0$. Then there exists an embedding of $X$ into $Y$ if and only if
\begin{align}
\# \{p \in P(X): \per(p) = n \}\leq \# \{p \in P(Y) : \per(p) = n\},
 \quad \quad \quad  n \in \Bbb N.
\end{align}
\end{Embedding Theorem}
\begin{proof}
We can assume that the irreducible and aperiodic shift  of finite type $Y$ is given as the edge shift of a finite strongly connected aperiodic directed graph 
$G(\frak V_Y, \mathcal E_Y)$.
Let
$$
p = (p_k)_{k \in \Bbb Z}
$$
be a periodic point of $Y$. Denote the smallest period of $p$ by $\pi$  and let 
$\sigma \in \Sigma$
such that $\sigma$ can follow $p_{\pi - 1}$ and is distinct from ${p_0}$. Let $M, K > 1$ and set
$$
\ell = MK\pi.
$$

Also set
$$
a = p_{[0, \pi)},
$$
and set for $\xi, \zeta \in \Sigma$,
\begin{align*}
\mathcal D_l( \xi, \zeta) = \{ (d_ r)_{1\leq k \leq q }\in
\Gamma^{\langle + \rangle}(\xi) \cup  \Gamma^{\langle - \rangle} (\zeta):
d_{[r. r +M\pi)} \not = a^M, 1 \leq r \leq q - M\pi\}, 
\\
 2\ell \leq q + M\pi +1 \leq 4 \ell.
\end{align*}

We require that
\begin{align}
\#(\mathcal L_l(X)) \leq \#( \mathcal D_{l- M \pi -1}( \xi, \zeta)  ), \quad \quad \quad
 2\ell \leq l + M\pi + 1 \leq 4\ell,  \ \ \xi, \zeta \in \Sigma.
\end{align} 

It is known from general considerations (see e.g. \cite{L}) that the parameters $M$
 and $K$ can be chosen such that (2) holds. For instance, let
 $$
 K \geq \frac{h\pi }{\Delta},
 $$
 and, with an $R \in \Bbb N$ such that one can connect any two vertices of 
 $G(\frak V_Y, \mathcal E_Y)$ by a path of length $R$,  let 
 $$
 M \geq \frac{ h(\pi + R)}{\Delta},
 $$
 such that
 $$
 \# (\mathcal L_q(X)))< e^{(h - 2\Delta)q}, \quad \quad \quad q >(M-1)\pi.
 $$
 
  Set
 $$
 \alpha = p_0.
 $$
 By  (2) there exist injections
 \begin{align}
&\Xi_l(\sigma, \alpha): \mathcal L_l(X) \to \mathcal D_{l - M\pi -1}(\sigma, \alpha), \quad \quad
2 \ell  < l + M\pi +1 \leq 4\ell,
\\
&\widetilde {\Xi}(\zeta): \mathcal L_\ell(X) \to \mathcal D_{\ell}(\zeta, \alpha), \quad \quad \zeta \in \Sigma.
\end{align}

By (1) one can choose an embedding
$$
\psi_\circ: (P(X), S_X\restriction P(X)) \hookrightarrow (P(Y), S_X\restriction P(Y)).
$$      
For such a choice assign for all $n\in \Bbb N$ to every orbit $\frak p$ of of length $n$ 
of $X$ 
 an orbit $\psi_\circ(\frak p) $ of length $n$ of $Y$, choose $p \in  \frak p $ and
$ \psi_\circ(p) \in \psi_\circ(\frak p) $,
 and set
 $$
 \psi_\circ(S_X^k p) = S_X^k (\psi_\circ( p) ),  \quad \quad 0 \leq k < n.
 $$
 
Set for $L> 1$
$$
\mathcal B_{\ell, L}(X) = \bigcup_{1 \leq l \leq \ell} \{p_{[0, \ell L)} : p \in   P(X), \per(p) = l\}.
$$
Denote for $b \in\mathcal B_{\ell, L}(X) $ by $p^{(b)}$ the periodic point such that
$$
b = p_{[0, \ell L)}.
$$
 The embedding $\psi_\circ$ induces an injection
 $$
 \Psi:  \mathcal B_{\ell, 2}(X)    \hookrightarrow  \mathcal B_{\ell, 2}(Y) 
 $$
 by
 \begin{align}
 \Psi(b)= \psi_\circ(p^{(b)})_{[0, 2\ell)}.
 \end{align}
 
  Choose $L> 1$ such that 
  the set $\mathcal B_{\ell, L}(X)$ contains every periodic word in $\mathcal L(X)$ of length 
  $\ell L$
 and period $l, 1\leq l \leq \ell$.
Set
 $$
 \mathcal C_{\ell, L}(X)  = \mathcal L_{\ell L}(X) \smallsetminus \mathcal B_{\ell, L}(X),
 $$
and enumerate
$$
 \mathcal C_{\ell, L}(X)  =\{c^{(k)}:1 \leq k \leq \#(\mathcal C_{\ell, L}(X)) \}.
$$
Set
$$
Z^{(k)} = \{c \in  X: x_{[0, \ell L))}  = c^{(k)}\}, \quad  1 \leq k \leq \#(\mathcal C_{\ell, L}(X)),
$$
and set
$$
F^{(1)} = Z{(1)},
$$
\begin{align}
F^{(k)} = Z^{(k)} \smallsetminus (\bigcup_{- \ell < l < \ell, 1 \leq m <k} S_X^l F{(k)}   ),
1< k \leq \#(\mathcal C_{\ell, L}(X)),
\end{align}
$$
F = \bigcup_{1\leq k \leq \#(\mathcal C_{\ell, L}(X))} F^{(k)} .
$$

Concerning the set $F$ there are two observations. 
The first observation: If $x \in X, k \in \Bbb N,$ are such that
\begin{align}
x \in F \cap S_X^k F,
\end{align}
then $k\geq \ell.$  Indeed, assume that 
\begin{align}
k < \ell.
\end{align}
By (6) one has $q, r$ such that
$$
x \in F^{(q)} \cap S_X^k F^{(r)}.
$$
Equality of $q$ and $r$ would imply that $x_{[0, \ell L)}$ is in $\mathcal B_{\ell  L}$, which by (7) is not the case. If $q\not = r,$ then
$$
 F^{(q)} \cap S_X^k F^{(r)} \not = \emptyset
$$
which contradicts (6) by (8).
The second observation: If $x \in X$ such that
\begin{align}
S_X^k x \not \in F,   \quad \quad - \ell < k < \ell,
\end{align}
then the word $x_{[0, \ell L)}$ is in $\mathcal B_{\ell  L}$. Indeed, assume that 
there is an $m$ such that
$$
x_{[0, \ell L)} = c^{(m)}.
$$
By (9) $x \in F^{(m)}$ is impossible. This would mean that   
$$
x \in \bigcup_{1 \leq r < m, -\ell <k <\ell} S_X^k F{(r)},
$$
which contradicts (6).

These observations having been made we use the injections (3 - 4) and (5) in the formulation of   coding instructions for an embedding $\psi: X \hookrightarrow Y$ that extends the 
embedding $\psi_\circ: P(X) \hookrightarrow P(Y)$. With the the notation
$$
y = \psi_\circ(x),  \quad  x \in  X,
$$
and with the notation  $\xi$ and $\zeta$ for the first and last symbols of 
$\Psi (x_{[0, 2 \ell)})$
the instructions are:

 For
 $$
 0 < l <3\ell
 $$ 
 and
  $$
  x \in S_X^\ell F \cap \bigcap_{-\ell < k < l} (X \smallsetminus S_X^{-k}F) \cap S_X^{-l} F
  $$
  set
  $$
  y_{[-\ell, l)} = (a^M\sigma, \Xi_{\ell + l}(\sigma,\alpha))( x_{[-\ell, l)})).
  $$
 
 For 
  $$
  x \in S_X^{-\ell} F \cap \bigcap_{-\ell < k < 3\ell} (X \smallsetminus S_X^{-k}F)
  \cap S_X^{3\ell} F
  $$
  set
  $$
  y_{[-\ell, 3\ell)} = 
  (a^M\sigma, \Xi_\ell(\sigma, \xi)( x_{[0, 2 \ell)}), \Psi (x_{[0, 2 \ell)}),
  \widetilde{\Xi}( \zeta , \alpha)(x_{[2\ell, 3 \ell)} )).
  $$

 For 
  $$
  x \in  \bigcap_{-\ell \leq l < 3\ell } (X \smallsetminus S_X^{-k}F) 
  $$
   set
  $$
  y_{[0, 2\ell)} = \Psi(x_{[0, 2\ell)}).
  $$
  
  For 
  $$
  x \in  \bigcap_{-\ell \leq k < 3l} (X \smallsetminus S_X^{-k}F)  \cap S_X^{-3\ell}
  $$
   set
  $$
  y_{[0, 3\ell)} = ( \Psi (x_{[0, 2 \ell)}),  \widetilde{\Xi}( \zeta) (x_{[2\ell, 3 \ell)})).
  $$
  
  The point $x$ can be reconstructed from the point $y$:
  The word $a^M\sigma $ does not overlap itself.
  By construction one has for $x \in X$ that $x \in F$ precisely if
  $$
  y _{[0, M \pi + 1)} =a^M\sigma.
  $$
This circumstance and the injectivity of the mappings in   and  are sufficient for the reconstruction of $x$ from $\psi (x)$.
\end{proof}
   
For the proof of Theorem 1 we require that
\begin{align}
\#(\frak V ) \ \#(\mathcal L_l(X)) \leq \#( \mathcal D_{l- M \pi 1}( \xi, \zeta)  ), \quad \quad \quad
 2\ell \leq l + M\pi + 1 \leq 4\ell
\end{align} 
The parameters $M$
 and $K$ can be chosen such that (010) holds. If. is satisfied then the injections.  can be chosen such that the projection is bijective on the image of $X$ under the embedding of $X$ into $Y$.

\par\noindent Wolfgang Krieger
\par\noindent Institute for  Mathematics, 
\par\noindent  University of Heidelberg,
\par\noindent Im Neuenheimer Feld 205, 
 \par\noindent 69120 Heidelberg,
 \par\noindent Germany
\par\noindent krieger@math.uni-heidelberg.de

\end{document}